\documentclass[a4paper,12pt]{amsart}
\usepackage{amsmath,amsfonts,amssymb,amsthm,latexsym,amscd}
\usepackage[dvips]{graphicx}
\usepackage[T1]{fontenc}

\usepackage{parskip}

\def\D{\mathbb{D}}
\def\Diff{\operatorname{Diff}}

\def\a{\alpha}
\def\diam{\operatorname{diam}}
\def\g{\gamma}
\def\o{\omega}

\def\vol{\operatorname{vol}}

\newtheorem{thm}{Theorem}[section]
\newtheorem{thm*}{Theorem}
\newtheorem{lem}[thm]{Lemma}

\newtheorem*{q*}{Question}
\newtheorem*{rem*}{Remark}

\theoremstyle{definition}
\newtheorem{ex}[thm]{Example}

\newtheorem{rem}[thm]{Remark}

\def\OP{\operatorname}
\def\B{\mathbf}

\begin{document}

\title[Quasi-isometric embeddings]{Quasi-isometric embeddings into diffeomorphism groups}

\author{Michael Brandenbursky and Jarek K\k{e}dra}

\begin{abstract}
Let $M$ be a smooth compact connected oriented manifold of dimension at least
two endowed with a volume form. Assuming certain conditions on the
fundamental group $\pi_1(M)$ we construct quasi-isometric embeddings of
either free Abelian or direct products of non-Abelian free groups into
the group of volume preserving diffeomorphisms of $M$ equipped with
the $L^p$ metric induced by a Riemannian metric on $M$.

\end{abstract}

\maketitle

\section{Introduction}

\subsection{The $L^p$-metric}
Let $M$ be a compact connected and oriented Riemannian manifold and
let $\Diff(M,\mu)$ denote the group of $C^r$--diffeo\-mor\-phisms of
$M$ acting by the identity on a neighborhood of the boundary and
preserving the volume form $\mu$ induced by the metric.
Unless otherwise stated we assume that $\Diff(M,\mu)$ is
equipped with the $C^{k}$--topology for some fixed
$1\leq k\leq r\leq \infty$.

In the present paper we study the geometry of the identity component
$\Diff_0(M,\mu)$ of the above group endowed with the right
invariant $L^p$-metric. It is defined as follows. Let
$$
\mathcal{L}_p\{g_t\}:=
\int_0^1 dt \left(\int_M|\dot{g_t}(x)|^p\mu \right)^{\frac 1p}
$$
be the $L^p$-length
of a smooth isotopy $\{g_t\}_{t\in [0,1]}\subset\Diff_0(M,\mu)$,
where $|\dot{g_t}(x)|$ denotes the length of the tangent
vector $\dot{g_t}(x)\in T_xM$ induced by the Riemannian
metric. Observe that this length is right-invariant, that is,
$\mathcal{L}_p\{g_t\circ f\}=\mathcal{L}_p\{g_t\}$
for any $f\in \Diff(M,\mu)$. It defines a right-invariant
metric on $\Diff_0(M,\mu)$ by
$$
{\bf d}_p(g_0,g_1):=\inf_{g_t}\mathcal{L}_{p}\{g_t\},
$$
where the infimum is taken over all paths from
$g_0$ to $g_1$.

If $p=2$ then the group $\Diff_0(M,\mu)$ is in fact equipped with a
Riemannian metric inducing the above $L^2$-length.  The geodesics of
this metric are the solutions of the equations of the flow of an
incompressible fluid \cite{Ar}, which makes the $p=2$ case the most
interesting. It is known that if $M$ is a simply connected Riemannian
manifold of dimension at least three then the $L^2$-diameter of the
group $\Diff_0(M,\mu)$ is finite \cite{Sh}. On the other hand
Eliashberg and Ratiu \cite{ER} proved that this diameter is infinite
for surfaces and for manifolds with positive first Betti number.  See
Arnol'd-Khesin \cite{AK} and Khesin-Wendt \cite[Section~3.6]{KW} for
a detailed discussion.

\begin{rem}
The property of preserving the volume is essential
to prove the right invariance of the above metric.
One can define the $L^p$-metric by defining first a norm
of a diffeomorphism $g$ by
$$
\|g\|_{p}:=\inf_{g_t}\mathcal{L}_{p}\{g_t\},
$$
where the infinum is taken over all isotopies from the
identity to $g$. Then the metric
${\bf d}_p(g,h):=\|gh^{-1}\|_p$ is right invariant by
definition. However, in order to prove the triangle
inequality it is necessary to use the property of
preserving the volume.
\end{rem}

\subsection{The main result}

%Before stating our main result we want to make the following important remark.

%\begin{rem}\label{R:vol}
%\hfill
%\begin{itemize}
%\item
%The volume preserving property is used to prove
%the triangle inequality and the right-invariance of the $L^p$-metric.
%The only place where we use that the group $\Diff(M,\mu)$ is the group of
%\emph{volume-preserving} diffeomorphisms is in the proof of right-invariance of
%$L^p$-metric on $\Diff_0(M,\mu)$. More precisely, in order to prove that
%${\bf d}_p(\Id,g)={\bf d}_p(h,gh)$ we need $h$ to be volume-preserving.
%\item
%Throughout the paper we are restricting ourselves to
%the group of smooth diffeomorphisms. However, the same proof
%works for the group of $C^k$ diffeomorphisms with
%$C^k$-topology for $k\geq 1$.
%We equip the group $\Diff(M,\mu)$ with the Whitney $C^{\infty}$-topology,
%because in this topology any smooth isotopy between two diffeomorphisms in
%$\Diff(M,\mu)$ is a continuous map (see \cite{Ba}), and we use smooth
%isotopies in order to define $L^p$-metric on~$\Diff_0(M,\mu)$.
%\end{itemize}
%\end{rem}

A map $\psi\colon (X_1,\OP{d}_1)\to (X_2,\OP{d}_2)$ between metric spaces
is called {\em large scale Lipschitz} \cite[Remark 1.9]{Roe}
if there exist constants
$A,B\geq 0$ such that
$$
\OP{d}_2(\psi(x),\psi(y))\leq A\cdot \OP{d}_1(x,y)+B.
$$

Let $m\in M$ be a reference point.
Let $\OP{ev}_m\colon \Diff_0(M,\mu)\to M$ be the evaluation map
defined by $\OP{ev}_m(f):=f(m)$ and let $G_{\mu}\subset \pi_1(M)$ be
the image of the homomorphism induced by $\OP{ev}_m$.  It is easy to
prove that $G_{\mu}$ is contained in the center of $\pi_1(M)$.  The
subgroup $G_{\mu}$ is called the {\em Gottlieb group} associated with
the volume form $\mu$ because the groups of similar origin were first
studied by Gottlieb in \cite{Go}.

Let $\Diff(M,\mu,m)\subset \Diff_0(M,\mu)$ be the isotropy of
the reference point $m$. Let us define a homomorphism
$$
\Phi\colon \Diff(M,\mu,m)\to \pi_1(M)/G_\mu
$$ as follows. Let $g\in \Diff(M,\mu,m)$ and let $\{g_t\}\subset
\Diff_0(M,\mu)$ be a smooth isotopy from the identity to $g$. The
value $\Phi(g)$ is represented by the loop $\{g_t(m)\}$ and it is
straightforward to show that $\Phi$ is well defined.

Let $B(m,r)\subset M$ be a ball of radius $r>0$ centered at
a reference point $m\in M$. Let $\Diff(M,\mu,B(m,r))$ denote the
subgroup of $\Diff_0(M,\mu)$ consisting of diffeomorphisms
preserving the ball $B(m,r)$ pointwise. The metric on the group $\Diff(M,\mu,B(m,r))$ is induced from the $L^p$-metric on $\Diff_0(M,\mu)$.
The following is our main technical result which is
proven in Section \ref{S:proofs}.

\begin{thm}\label{T:lsl}
Let $M$ be a compact connected and oriented Riemannian manifold.
For all small enough $r>0$ the homomorphism
$$
\Phi\colon \Diff(M,\mu,B(m,r)) \to \pi_1(M)/G_\mu
$$
is surjective and large scale Lipschitz with respect to the
$L^p$-metric on the group $\Diff(M,\mu,B(m,r))$ and the
word metric on $\pi_1(M)/G_{\mu}$.
\end{thm}

Recall that the word norm on a group $\Gamma$ generated
by a symmetric finite set $S\subset \Gamma$ is defined by
$$
|\g|_S:=
\min\{k\in \mathbb N\,|\,\g=s_1\ldots s_k \text{ where }s_i\in S\}.
$$
The word metric is defined by $\OP{d}_S(\g_1,\g_2):=|\g_1(\g_2)^{-1}|_S$.
It is right-invariant and it depends on the choice of a finite
generating set up to a quasi-isometry (defined below)
\cite[Example 8.17]{BH}.

\subsection{Applications}\label{SS:applications}

Let $(X_1,\OP{d}_1)$ and $(X_2,\OP{d}_2)$ be two metric spaces. A function
$f\colon X_1\to X_2$
is a \emph{quasi-isometric embedding} if there exist
two constants $A\geq1$ and $B\geq0$ such that
\begin{equation}\label{eq:q-i-e}
\frac{1}{A}\OP{d}_1(x,y)-B\leq \OP{d}_2(f(x),f(y))\leq A\OP{d}_1(x,y)+B.
\end{equation}
In the case when $(X_1,\OP{d}_1)$ and $(X_2,\OP{d}_2)$ are metric groups, we
require $f$ to be an \emph{injective homomorphism}.  We say that $f$ is a
\emph{quasi-isometry} if, in addition to \eqref{eq:q-i-e},
there exists a constant $C\geq0$ such that for every $u\in X_2$ there
exists $x\in X_1$ with the property
$$\OP{d}_2(u,f(x))\leq C.$$

Let $\Gamma$ be a group equipped with
the word metric associated with a finite generating set $S\subset \Gamma$.
An element $\g$ is called {\em undistorted}
in $\Gamma$ if there exists a positive constant $C>0$ such that
$$
|\g^n|_S\geq C\cdot n.
$$
Otherwise, $\g$ is called {\em distorted}.
These properties do not depend on the choice of a finite
generating set.

\begin{thm}\label{T:undistorted}
Let $M$ be a compact connected and oriented Riemannian manifold of
dimension at least two. If $\pi_1(M)/G_{\mu}$ contains an undistorted
element then $(\Diff_0(M,\mu),{\bf d}_p)$ contains quasi-isometrically
embedded free Abelian group of an arbitrary rank. In particular,
the metric group $(\Diff_0(M,\mu),{\bf d}_p)$ has infinite diameter.
\end{thm}

This result generalizes a theorem of Eliashberg and Ratiu \cite{ER}
where they prove the infiniteness of the diameter under the assumption
that the first Betti number of $M$ is positive and the center of the
fundamental group is trivial. Notice that there exist compact
oriented manifolds with the fundamental group isomorphic to an
arbitrary finitely presented group. Such a group can be chosen to have
finite abelianization (hence the first Betti number is equal to zero)
and usually groups have undistorted elements. For example,
$\OP{SL}(2,\B Z)$ has finite abelianization.
Moreover, if $g\in \OP{SL}(2,\B Z)$ has an eigenvalue
$\lambda$ such that $|\lambda|\neq 1$ then the cyclic group generated
by $g$ is undistorted.

In fact, no finitely presented group
with all its elements distorted is known.
On the other hand, Osin \cite{Os} constructs an
example of an infinite finitely generated group with exactly two
conjugacy classes.  In particular, such a group does not have
undistorted elements.

\begin{thm}\label{T:free}
Let $M$ be a compact connected and oriented Riemannian manifold of
dimension at least three. If $\pi_1(M)/G_{\mu}$ contains
quasi-isometrically embedded non-Abelian free group then
$(\Diff_0(M,\mu),\B d_p)$ contains quasi-isometrically
embedded direct product of any finite number of
free groups of arbitrary ranks.
\end{thm}

Essentially, the idea of the proof is to appropriately embed
disjoint copies of the figure eight into $M$ and suitably
apply Theorem \ref{T:lsl}. Thus the situation is a bit
different for surfaces, where we get the following slightly
weaker statement.

\begin{thm}\label{T:free_surface}
Let $\Sigma_{g,k}$ be a compact connected and oriented
surface of genus $g$ with $k$ boundary components.
Then
$(\Diff_0(\Sigma_{g,k},\mu),\B d_p)$ contains quasi-isometrically
embedded direct product of $2g+k-2$ copies of
finitely generated non-Abelian free groups of arbitrary ranks.
\end{thm}

\begin{rem}
Let $\D^2$ be a unit Euclidean disc in $\B R^2$, i.e. in our notation
$\D^2$ is diffeomorphic to $\Sigma_{0,1}$. In \cite{BG} Benaim and
Gambaudo showed, using a different method, that the group
$(\Diff(\D^2,\mu),\B d_p)$ contains quasi-isometrically embedded
finitely generated free or free Abelian group of arbitrary rank. Crisp
and Wiest proved the same fact for planar right-angled Artin groups
\cite{CW}. We also would like to mention that M. Kapovich proved that
any right-angled Artin group embeds into the group of Hamiltonian
diffeomorphisms of any symplectic manifold $(M,\o)$, see
\cite{Ka}. However, it is not known whether the embedding he
constructs in \cite{Ka} is quasi-isometric with respect to
$L^p$-metric.
\end{rem}

\subsection{Examples}\label{SS:examples}

\begin{ex}\label{E:p3}
Since the Artin pure braid group $\B P_3$ on three strands is
isomorphic to $\B F_2\times {\B Z}$ (see the proof of
Theorem 1.16 in \cite{KT}) it embeds quasi-isometrically
into $(\Diff_0(M,\mu),\B d_p)$, where the manifold $M$ is as in Theorem \ref{T:free}.
\end{ex}

\begin{ex}\label{E:stallings}
If a subgroup $\Gamma$ of a direct product of $n$
free groups has finitely generated homology up to
degree $n$ then $\Gamma$ contains a finite index
subgroup isomorphic to a direct product of at most
$n$ free groups~\cite{BHMS}.
On the other hand, there are examples of finitely
presented subgroups of a direct product of free
groups which have more complicated finiteness
properties. More concretely, the kernel $\Gamma$
of the homomorphism
$(\B F_2)^n \to \B Z$
sending each generator to one is of type
$\OP{FP}_n$ but not $\OP{FP}_{n+1}$.
These examples are known as {\em Bieri-Stallings groups},
see Bestvina-Brady~\cite[Example 6.3]{BB}.

Suppose that  $\Gamma \to \B F_2\times \ldots \times \B F_2$
is the inclusion of a finitely presented Bieri-Stallings group.
It follows from the proof of Theorem 11.7 of Dison's
thesis \cite{Di} that this inclusion is a quasi-isometric
embedding. Indeed, since the quotient is isomorphic to
$\B Z$, its isoperimetric function is linear and hence,
by Dison's Lemma 9.5 in \cite{Di} the distortion function of the
above inclusion is linear. A more straightforward proof can
be found in Bridson-Haefliger \cite[Exercise 5.12(3)]{BH}.
\end{ex}

Recall that if $\psi\colon H\to G$ is an injective homomorphism
of metric groups then its distortion function
$\Delta\colon \B R_+\to \B R_+$ is defined by
$$
\Delta(r) := \max\{\OP{d}_H(1,h)\,|\, \OP{d}_G(1,\psi(h))<r\}.
$$
Observe that the distortion function is linear if and only if
$\psi$ is a quasi-isometric embedding.

\begin{rem}\label{rem:torus}
Let $M=\B T^n$ be the $n$-dimensional torus equipped with
a volume form $\mu$.
Neither the result of Eliashberg and Ratiu nor
our Theorem \ref{T:undistorted} apply to $\Diff(\B T^n,\mu)$ because
$G_{\mu}=\pi_1(\B T^n)$ and the whole group is central.
\end{rem}

\begin{ex}\label{E:domain}
Let $M\subset \B R^3$ be a closed domain with free non-Abelian
fundamental group. The geometry of $(\Diff_0(M,\mu),\B d_p)$
models the behavior of an incompressible fluid
filling a tank of the shape of~$M$. Theorem \ref{T:free}
describes a large scale complexity of mixing such a fluid.
\end{ex}

\subsection{The symplectic case}\label{SS:symp}

If $(M,\omega)$ is a symplectic manifold, then
the group $\Diff_0(M,\mu)$ in all the results above
can be replaced either by the group $\OP{Symp}_0(M,\omega)$
of symplectic diffeomorphisms isotopic to the identity,
or by the group $\OP{Ham}(M,\omega)$
of Hamiltonian diffeomorphisms, see Remark \ref{R:symp}
in the proof of Theorem \ref{T:lsl}.

\section{Proofs}\label{S:proofs}

\subsection{An abstract lemma}\label{SS:abstract}

Let $(G,\OP{d}_G)$ be a metric group with the identity element
$1_G$. For $g\in G$ we set
$$\|g\|_G:=\OP{d}_G(1_G,g)\quad\textrm{and}\quad \diam(G):=\sup\limits_{g\in G}\|g\|_G.$$

\begin{lem}\label{L:abstract}
Let $(G,\OP{d}_G)$, $(H,\OP{d}_H)$, $(K,\OP{d}_K)$ be three metric groups, such
that $H$ is finitely generated and $\OP{d}_H$ is a word metric w.r.t. some
finite generating set $S$. Suppose that $\Phi\colon G\to K$ is a large scale
Lip\-schitz homomorphism. Let $\Psi\colon H\to G$ be a homomorphism, such that
$\Phi\Psi\colon H\to K$ is a quasi-isometric embedding.  Then $\Psi$ is a
quasi-isometric embedding.
\end{lem}

\begin{proof}
The homomorphism $\Psi$ is injective, because $\Phi\Psi$ is
injective. Let $h\in H$. The homomorphism $\Phi\Psi$ is a
quasi-isometric embedding, hence there exist two constants $A_1\geq0$
and $B_1\geq0$ such that

\begin{equation}\label{eq:ab-lem-1}
A_1\|h\|_H-B_1\leq\|\Phi\Psi(h)\|_K
\end{equation}

The homomorphism $\Phi$ is a large scale Lipschitz, hence there exist
two constants $A_2>0$ and $B_2>0$ such that

\begin{equation}\label{eq:ab-lem-2}
\|\Phi\Psi(h)\|_K\leq A_2\|\Psi(h)\|_G+B_2.
\end{equation}

Combining inequalities \eqref{eq:ab-lem-1} and \eqref{eq:ab-lem-2} we get

\begin{equation}\label{eq:ab-lem-3}
\frac{A_1}{A_2}\|h\|_H-\left(\frac{B_1+B_2}{A_2}\right)\leq\|\Psi(h)\|_G.
\end{equation}

Let $k$ be such that $\|h\|_H=k$. The group $H$ is finitely generated,
hence $h=h_1\cdot\ldots\cdot h_k$, where $h_i\in S$ for each
$1\leq i\leq k$.  Denote by
$$M_\Psi:=\max\{\|\Psi(h_1)\|_G,\ldots,\|\Psi(h_k)\|_G\}.$$
It follows that
\begin{equation}\label{eq:ab-lem-4}
\|\Psi(h)\|_G=
\|\Psi(h_1)\cdot\ldots\cdot\Psi(h_1)\|_G\leq\sum\limits_{i=1}^k\|\Psi(h_i)\|_G\leq M_\Psi\cdot k
=M_\Psi\|h\|_H.
\end{equation}
Inequalities \eqref{eq:ab-lem-3} and \eqref{eq:ab-lem-4} conclude the proof of the lemma.
\end{proof}

In the proofs below we use the fact that the metric on the groups $\Diff(M,\mu, B(m,r))$ and $\Diff(M,\mu, \sqcup_i B(m_i,r))$
(this group is defined in the proof of Theorem \ref{T:undistorted}) is induced from the $L^p$-metric on the whole group $\Diff(M,\mu)$.

\subsection{Proof of Theorem \ref{T:lsl}}\label{SS:proof_lsl}
In the first part we prove the surjectivity of the homomorphism
$\Phi\colon \Diff(M,\mu,B(m,r)) \to \pi_1(M)/G_\mu$.

Let $S':=\{[\g_1],\ldots,[\g_k]\}\subset \pi_1(M)$ be a symmetric
generating set of the fundamental group of $M$ such that each
representative $\g_i$ is a simple closed curve. It follows from the
tubular neighborhood theorem that for each $\g_i$ there exists
$r_i>0$, the standard $(n-1)$-dimensional ball
$B_{r_i}^{n-1} \subset \B{R}^n$ of radius $r_i>0$,
and a volume-preserving embedding
$$
\OP{emb}_i:B_{r_i}^{n-1}\times \B S^1\hookrightarrow M,
$$
such that $\OP{emb}_i|_{\{0\}\times \B S^1}$ is the curve $\g_i$ and
$r_i=\OP{d}_M(\g_i, \OP{emb}_i|_{\partial B_{r_i}^{n-1}\times \B S^1})$.
The volume form on the product $B_{r_1}^{n-1}\times \B S^1$ is the
standard Euclidean volume and $\OP{d}_M$ denotes the distance on
$M$ induced by the Riemannian metric.

Let $r=\min\limits_{1\leq i\leq k}{r_i}$. Then
$\OP{emb}_i:B_r^{n-1}\times \B S^1\hookrightarrow M$
is volume-preserving for each $i$ and
$\OP{emb}_i|_{{0}\times  \B S^1}=\g_i$.
It is straightforward to
construct a smooth isotopy of volume-preserving diffeomorphisms
$$
g_t:B_{r}^{n-1}\times \B S^1\to B_{r}^{n-1}\times \B S^1
$$
between $g_0=\OP{Id}$ and $g_1$ such that:
\begin{itemize}
\item
For each $t\in [0,1]$ the diffeomorphism $g_t$ equals to the identity
in the neighborhood of $\partial B_{r}^{n-1}\times \B S^1$, and
the time-one map $g_1$ is equal to the identity on $B_{r'}^{n-1}\times \B S^1$,
where $0<r'<r$.
\item
Each diffeomorphism $g_t$ preserves the foliation of
$B_r^{n-1}\times \B S^1$ by the circles $\{x\}\times \B S^1$ and for
every $x\in B^{n-1}_{r'}$ the restriction
$g_t\colon \{x\}\times \B S^1\to \{x\}\times \B S^1$ is the rotation
by $2\pi\,t$. It follows that the time-one map $g_1$ equals the
identity on $B^{n-1}_{r'}\times \B S^1$.
\end{itemize}

We identify $B^{n-1}_{r}\times \B S^1$ with its image
with respect to the embedding $\OP{emb}_i$. Then we extend
each isotopy $g_t$ by the identity on
$M\setminus B_r^{n-1}\times \B S^1$
obtaining smooth isotopies $g_{t,i}\in \Diff_0(M,\mu)$.
This shows that every representative $\g_i$ of a generator
of the fundamental group of $M$ arises as simple closed
curve $\{g_{t,i}(m)\}$ and hence the homomorphism
$$
\Phi\colon \Diff(M,\mu,B(m,r))\to \pi_1(M)/G_{\mu}
$$
is surjective.

\begin{rem}\label{R:symp}
Notice that if $M$ is a symplectic manifold then the
above isotopies can be constructed to be Hamiltonian.
\end{rem}

Let $\Pi:\pi_1(M)\to\pi_1(M)/G_\mu$ be the projection homomorphism.
Consequently $S:=\Pi(S')$ is a finite generating set for the quotient
$\pi_1(M)/G_{\mu}$.  Let $\Pi_M\colon {M}_{\bullet}\to M$ be the
Riemannian covering associated with $\Pi$. This means that the metric
on $M_{\bullet}$ is induced from the Riemannian metric on $M$. The
corresponding distance will be denoted by $\OP{d}_{\bullet}$.

Now we shall prove that $\Phi$ is a large scale Lipschitz map.
That is, we show that there exist positive constants $A$ and $B$
independent of $g$ such that
$$
A\cdot \|g\|_p+B\geq\|\Phi(g)\|_S,
$$
where $\|g\|_p:={\bf d}_p(g,\OP{Id})$ is the $L^p$-norm
of the diffeomorphism $g$.

Let $g\in \Diff(M,\mu,B(m,r))$ and let $\{g_t\}_{t\in [0,1]}\in \Diff_0(M,\mu)$
be an isotopy from the identity to $g$. It follows from the H\"older
inequality that $\|g\|_p\geq C_p\cdot \|g\|_1$,
where $C_p$ is some positive constant independent of $g$. Hence it is
enough to prove the statement for $p=1$.

Let ${m}_{\bullet}\in \Pi_M^{-1}(m)$, and
let $\{g_{\bullet,t}(m_{\bullet})\}$ be the lift of
$\{g_t(m)\}$ starting at the point $m_{\bullet}$.
The manifold $M$ is compact, hence by the $\check{\textrm{S}}$varc-Milnor
lemma \cite{BH,Mil}, the inclusion of the orbit of $m_{\bullet}$
with respect to the deck transformation group $\pi_1(M)/G_{\mu}$ defines
a quasi-isometry
\begin{equation*}
\pi_1(M)/G_{\mu} \stackrel{q.i.}\simeq ({M}_{\bullet},\OP{d}_{\bullet}).
\end{equation*}
In particular, it means that there exist positive constants $A',B'$, such that
\begin{equation}\label{eq:q-i-Milnor}
{d}_{\bullet}({m}_{\bullet},{g}_{\bullet,1}({m}_{\bullet}))\geq A'\|\Phi(g)\|_S-B'.
\end{equation}

Let $x\in B(m,r)\subset M$. We claim that the length of the flow-line
$g_t(x)$ is bounded by the distance
$\OP{d}_{\bullet}(m_{\bullet},g_{\bullet,1}(m_{\bullet}))$ up to the
diameter of the ball. To see this, consider the lift of each flow-line
$g_t(x)$ starting at a ball of radius $r$ in $M_{\bullet}$ centered at
$m_{\bullet}$ and observe that such a lift ends in a ball of radius
$r$ centered at $g_{\bullet,1}(m_{\bullet})$.
Indeed, let
$\a\colon [0,1]\to B(m,r)$ be a path between $x$ and $m$. Then the
map
$$
H:[0,1]\times[0,1]\to M,
$$
defined by $H(t,s)=g_t(\a(s))$ is
a homotopy from $\{g_t(m)\}$ to $\{g_t(x)\}$.
Lifting this homotopy shows that the lift of
$\{g_t(x)\}$ finishes at the ball $B(g_{\bullet,1}(m_{\bullet}),r)$.
Finally, we obtain that
\begin{equation}\label{eq:L-min}
\OP{Length}(g_t(x)):=\int_0^1|\dot{g}_t(x)|dt
\geq \OP{d}_{\bullet}({m}_{\bullet},{g}_{\bullet,1}({m}_{\bullet}))-2r.
\end{equation}
as claimed.
Combining inequalities \eqref{eq:q-i-Milnor} and \eqref{eq:L-min} we get that
\begin{equation*}
\OP{Length}(g_t(x))\geq A'\|\Phi(g)\|_S-(B'+2r)
\end{equation*}
for every $x\in B(m,r)$.
Hence by Fubini theorem and the above inequality we have
\begin{eqnarray*}
\mathcal{L}_1(\{g_t\})&=&\int_0^1 dt\left(\int_{M}|\dot{g}_t(x)|\mu\right) \\
&=&\int_{M}\mu \left(\int_0^1|\dot{g}_t(x)|dt\right)\\
&\geq& \vol(B(m,r))\cdot \min_{x\in B(m,r)}\OP{Length}(g_t(x))\\
&\geq&  \vol(B(m,r))\cdot A'\|\Phi(g)\|_S-\vol(B(m,r))\cdot(B'+2r).
\end{eqnarray*}
Since the above inequalities hold for any isotopy
$\{g_t\}_{t\in[0,1]}$ between the identity and $g$,
we obtain that
\begin{equation*}
\|\Phi(g)\|_S\leq A\cdot \|g\|_p+B,
\end{equation*}
where $A=(C_p\cdot \vol(B(m,r))^{-1}\cdot A'$ and $B=\frac{B'+2r}{A'}$
and this concludes the proof.
\qed

\subsection{Proof of Theorem \ref{T:undistorted}}
\label{SS:proof_undistorted}

Let $n\in \B N$ be a positive integer.
Recall that we need to prove that there exists a
quasi-isometric embedding of a free Abelian group
of rank $n$ into $\Diff_0(M,\mu)$.

Assume first that the dimension of $M$ is at least three.
Let $m_1,\ldots,m_n$ be distinct points in the interior of $M$ and
let $r>0$ be such that the balls $B(m_i,r)$
of radius $r$ centered at $m_i$ are pairwise disjoint.
Let $\g_{i,j}$ be simple closed curves representing
the generators of $\pi_1(M,m_i)$. We also assume
that $\g_{i_1,j_1}$ is disjoint from $\g_{i_2,j_2}$
whenever $i_1\neq i_2$. We choose $r$ small enough
such that the tubular neighborhood of radius $r$
of the above generators are disjoint.

Let $G_i\subset \pi_1(M,m_i)$ be the corresponding Gottlieb
group. The groups $\pi_1(M,m_i)/G_i$ are pairwise isomorphic.
Let $\g_i \in \pi_1(M,m_i)/G_i$ be an undistorted element
which exists according to the hypothesis. Let
$$
h\colon \B Z^n\to \pi_1(M,m_1)/G_1\times \ldots \times \pi_1(M,m_n)/G_n
$$
be a homomorphism defined by
$$
h(k_1,\ldots,k_n):= (\g_1^{k_1},\ldots,\g_n^{k_n}).
$$
It immediately follows from the fact that each $\g_i$ is undistorted
that $h$ is a quasi-isometric embedding.
Let $\Diff(M,\mu,\sqcup_i B(m_i,r))\subset \Diff_0(M,\mu)$ be the subgroup consisting
of diffeomorphisms preserving the disjoint union of balls
$B(m_i,r)$ pointwise. Let
$$
\Phi_i\colon \Diff(M,\mu,B(m_i,r))\to \pi_1(M,m_i)/G_i
$$
be the homomorphism defined in Theorem \ref{T:lsl}.
Consider a homomorphism
$$
\widetilde \Phi\colon \Diff(M,\mu,\sqcup_i B(m_i,r))\to
\pi_1(M,m_1)/G_1\times \ldots \times \pi_1(M,m_n)/G_n
$$
which is the composition of the (diagonal) inclusion
$$
\iota\colon \Diff(M,\mu,\sqcup_i B(m_i,r))\hookrightarrow
\prod_i \Diff(M,\mu,B(m_i,r))
$$
followed by the product homomorphism
$$
\prod_i \Phi_i\colon \prod_i \Diff(M,\mu,B(m_i,r))\to \prod_i \pi_1(M,m_i)/G_i.
$$
Since the inclusion $\iota$ is an isometric embedding and the
$\prod_i\Phi_i$ is large scale Lipschitz, according to Theorem \ref{T:lsl},
we obtain that $\widetilde \Phi$ is a large scale Lipschitz homomorphism.

Let $g_i\in \Diff(M,\mu,B(m_i,r))$ be an element such that
$\Phi_i(g_i)=\g_i$ and $g_i$ is supported in the
union of the tubular neighborhoods of the loops
representing the generators of $\pi_1(M,m_i)$ constructed
in the beginning of the proof. It follows that
the supports of $g_i$ and $g_j$ are disjoint
if $i\neq j$. The existence of $g_i$ follows from the
proof of Theorem~\ref{T:lsl}.
Let
$$
\Psi\colon \B Z^n\to \Diff(M,\mu,\sqcup_iB(m_i,r))\subset \Diff_0(M,\mu)
$$
be defined by
$$
\Psi(k_1,\ldots,k_n):= g_1^{k_1}\circ \dots \circ g_n^{k_n}.
$$
It is well defined because $g_i$ have pair-wise disjoint supports.
Recall that we have that
$$
h=\widetilde \Phi\circ \Psi\colon \B Z^n \to \prod_i\pi_1(M,m_i)/G_i
$$
and we know that $h$ is a quasi-isometric embedding and
$\widetilde \Phi$ is large scale Lipschitz.
Consequently the map $\Psi$ is a quasi-isometric
embedding according to Lemma \ref{L:abstract}.

Let us now consider the two-dimensional case.  Let $M=\Sigma_{g,k}$ be
a compact oriented surface of genus $g$ with $k$ boundary
components. Observe that $\pi_1(\Sigma_{g,k})/G_{\mu}$ is trivial if
either $g=0$ and $k\leq 2$ or $g=1$ and $k=0$. Otherwise it is either
free non-Abelian group or the fundamental group of a closed oriented
surface. In each case it is straightforward to define an embedding
$$
\OP{emb}\colon \B S^1\times [0,2n] \to M
$$
such that each loop $\OP{emb}(\B S^1\times\{t\})$ represents
an undistorted element in $\pi_1(M,\OP{emb}(1,t))$.
Let $g_i\colon M\to M$, for $i=1,\ldots,n$
be an area preserving diffeomorphism satisfying each of
the following conditions:
\begin{itemize}
\item
it is supported in $\OP{emb}(\B S^1\times (2i-2,2i))$;
\item
it  preserves the ball $\OP{emb}(B_i)$, where
$B_i\subset \B S^1\times [0,2n]$ is a ball of
diameter one centered at $(1,2i-1)$;
\item
it is the time one map of an isotopy from the
identity which acts as the full rotation
on the loop $\OP{emb}(\B S^1\times \{2i-1\})$.
\end{itemize}

As in the previous part the homomorphism
$\Psi\colon \B Z^n\to \Diff_0(M,\mu)$ defined
by $\Psi(k_1,\ldots,k_n):=g_1^{k_1}\circ\dots\circ g_n^{k_n}$
is the required quasi-isometric embedding.
\qed

\subsection{Proof of Theorem \ref{T:free}} \label{SS:free}

This proof is a modification of the proof of Theorem
\ref{T:undistorted} for three dimensional $M$ where the cyclic group
$\B Z$ is replaced by a non-Abelian free group $\B F_2$ on two
generators.
More precisely, let $f_{i},g_{i}\in \Diff_0(M,\mu)$ be
diffeomorphisms satisfying each of the following conditions
(we use here the notation of the proof of Theorem~\ref{T:undistorted}):
\begin{itemize}
\item
the support of $f_{i}$ and $g_{i}$ is contained in the neighborhood
of the union of the loops $\g_{i,j}$;
\item
the images $\Phi_i(f_{i})$ and $\Phi_i(g_{i})$ generate
the free non-Abelian group in $\pi_1(M,m_i)/G_i$.
\end{itemize}
Such diffeomorphisms can be constructed in a similar way as
$g_i$'s in the proof of Theorem \ref{T:undistorted}.

Let $w\in \B F_2$ be a reduced word and given two
elements $f,g\in \Diff(M,\mu)$ let $w(f,g)$ denote
the induced diffeomorphism of $M$.
Let
$$
\Psi\colon \B F_2\times \dots \times \B F_2\to
\Diff(M,\mu,\sqcup_iB(m_i,r))\subset \Diff_0(M,\mu)
$$
be defined by
$
\Psi(w_1,\ldots,w_n):= w_1(f_{1},g_{1})\circ \dots \circ
w_n(f_{n},g_{n})
$.
As before $\Psi$ is a quasi-isometric embedding of
a product of free groups on two generators into
$\Diff_0(M,\mu)$.
Since $\B F_2$ contains quasi-isometrically embedded
a non-Abelian free group of an arbitrary finite
rank \cite{Harpe} the proof is finished.
\qed

\subsection{Proof of Theorem \ref{T:free_surface}} \label{SS:proof_free_surface}

The proof of the two dimensional case
of Theorem \ref{T:undistorted} amounts to constructing
a number of disjoint simple closed curves representing
an undistorted element in the fundamental group of $M$.
The present proof is analogous in the sense that
we need to construct an embedding of the disjoint union
of $2g+k-2$ copies of the figure eight into $M$ such that
each embedding induces a quasi-isometric embedding
$\B F_2\to \pi_1(M,m_i)$ for $i=1,\ldots,2g+k-2$.
We leave this straightforward construction as an
exercise to the reader.

The rest of the proof is similar to the other proofs.
That is, we construct relevant diffeomorphisms $f_i,g_i$
and observe that the map
$$
\Psi\colon \B F_2\times \dots \times \B F_2\to
\Diff(\Sigma_{g,k},\mu,\sqcup_iB(m_i,r))
\subset \Diff_0(\Sigma_{g,k},\mu)
$$
defined by $\Psi(w_1,\ldots,w_n):= w_1(f_{1},g_{1})\circ \dots \circ w_n(f_{n},g_{n})$
is a quasi-isometric embedding.
\qed

\subsection*{Acknowledgments}
We would like to thank Jim Howie and Mark Sapir for answering our
questions. Example \ref{E:stallings} is due to Jim Howie.
We also thank Martin Bridson for useful historical comments concerning
this example.

This work has been done during the first author stay
in Aberdeen. His visit was supported by the Coleman-Cohen foundation.

\bibliographystyle{alpha}

\bigskip

Department of Mathematics, Vanderbilt University, Nashville, TN 37240\\
\emph{E-mail address:} \verb"michael.brandenbursky@vanderbilt.edu"

\bigskip

Department of Mathematics, Aberdeen University, Aberdeen, UK\\
\emph{E-mail address:} \verb"kedra@abdn.ac.uk"

\end{document}